\documentclass[12pt]{article}
\usepackage[a4paper,margin=2.5truecm]{geometry}
\usepackage{amsmath,amsfonts,amssymb,amsthm,txfonts,hyperref,url,tikz,graphicx}
\tikzset{>=stealth}
\newtheorem{lemma}{Lemma}
\newtheorem{conjecture}{Conjecture}
\newtheorem{problem}{Open Problem}
\newtheorem{observation}{Observation}
\newtheorem{theorem}{Theorem}
\newcommand{\N}{\mathbb N}
\newcommand{\Z}{\mathbb Z}
\newcommand{\Q}{\mathbb Q}
\title{A hypothetical upper bound on the heights of the solutions of a Diophantine equation with a finite number of solutions}
\author{Apoloniusz Tyszka}
\begin{document}
\begin{sloppypar}
\date{}
\maketitle
\begin{abstract}
Let \mbox{$f(1)=1$}, and let \mbox{$f(n+1)=2^{\textstyle 2^{\textstyle f(n)}}$} for every positive integer $n$.
We conjecture that if a system \mbox{${\cal S} \subseteq \{x_i \cdot x_j=x_k:~i,j,k \in \{1,\ldots,n\}\} \cup
\{x_i+1=x_k:~i,k \in \{1,\ldots,n\}\}$} has only finitely many solutions in \mbox{non-negative} integers
\mbox{$x_1,\ldots,x_n$}, then each such solution \mbox{$(x_1,\ldots,x_n)$} satisfies \mbox{$x_1,\ldots,x_n \leqslant f(2n)$}.
We prove: {\tt (1)} the conjecture implies that there exists an algorithm which takes as input a Diophantine equation,
returns an integer, and this integer is greater than the heights of integer (\mbox{non-negative} integer,
positive integer, rational) solutions, if the solution set is finite;
{\tt (2)} the conjecture implies that the question whether or not a Diophantine
equation has only finitely many rational solutions is decidable with an oracle
for deciding whether or not a Diophantine equation has a rational solution;
{\tt (3)} the conjecture implies that the question whether or not a Diophantine
equation has only finitely many integer solutions is decidable with an oracle
for deciding whether or not a Diophantine equation has an integer solution;
{\tt (4)} the conjecture implies that if a set \mbox{${\cal M} \subseteq \N$} has a \mbox{finite-fold}
Diophantine representation, then ${\cal M}$ is computable.
\end{abstract}
\mbox{finite-fold} Diophantine representation,
oracle for deciding whether or not a Diophantine equation has an integer solution,
oracle for deciding whether or not a Diophantine equation has a rational solution.
\section{Introduction and basic lemmas}
The height of a rational number \mbox{$\frac{p}{q}$} is denoted by \mbox{$h\left(\frac{p}{q}\right)$}
and equals \mbox{$\max(|p|,|q|)$} provided \mbox{$\frac{p}{q}$} is written in lowest terms.
The height of a rational tuple \mbox{$(x_1,\ldots,x_n)$} is denoted by \mbox{$h(x_1,\ldots,x_n)$}
and equals \mbox{$\max(h(x_1),\ldots,h(x_n))$}.
In this article, we present a conjecture which positively solves the following two open problems:
\begin{problem}\label{pro1}
Is there an algorithm which takes as input a Diophantine equation, returns an integer,
and this integer is greater than the moduli of integer (\mbox{non-negative} integer, positive integer)
solutions, if the solution set is finite?
\end{problem}
\begin{problem}\label{pro2}
Is there an algorithm which takes as input a Diophantine equation, returns an integer,
and this integer is greater than the heights of rational solutions, if the solution set is finite?
\end{problem}
\begin{lemma}\label{lem1}
For every \mbox{non-negative} integers $b$ and $c$, \mbox{$b+1=c$} if and only if
\mbox{$2^{\textstyle 2^b} \cdot 2^{\textstyle 2^b}=2^{\textstyle 2^c}$}.
\end{lemma}
\par
Let $[\cdot]$ denote the integer part function.
\begin{lemma}\label{lem2}
For every \mbox{non-negative} real numbers $b$ and $c$, \mbox{$b+1=c$} implies that
\mbox{$2^{\textstyle 2^{[b]}} \cdot 2^{\textstyle 2^{[b]}}=2^{\textstyle 2^{[c]}}$}.
\end{lemma}
\begin{proof}
For every \mbox{non-negative} real numbers $b$ and $c$, \mbox{$b+1=c$} implies that \mbox{$[b]+1=[c]$}.
\end{proof}
\section{A conjecture on the arithmetic of \mbox{non-negative} integers (rationals)}
Let 
\[
G_n=\{x_i \cdot x_j=x_k:~i,j,k \in \{1,\ldots,n\}\} \cup \{x_i+1=x_k:~i,k \in \{1,\ldots,n\}\}
\]
Let \mbox{$f(1)=1$}, and let \mbox{$f(n+1)=2^{\textstyle 2^{\textstyle f(n)}}$} for every positive integer $n$.
Let \mbox{$\theta(1)=0$}, and let \mbox{$\theta(n+1)=2^{\textstyle 2^{\textstyle \theta(n)}}$} for every positive integer $n$.
\begin{conjecture}\label{con1}
If a system \mbox{${\cal S} \subseteq G_n$}
has only finitely many solutions in \mbox{non-negative} integers \mbox{$x_1,\ldots,x_n$},
then each such solution \mbox{$(x_1,\ldots,x_n)$} satisfies \mbox{$x_1,\ldots,x_n \leqslant f(2n)$}.
\end{conjecture}
\par
An analogous conjecture in \cite[p.~710]{Tyszka1} has a much smaller bound.
Observations~\ref{obs1} and \ref{obs2} justify Conjecture~\ref{con1}.
\begin{observation}\label{obs1}
For every system \mbox{${\cal S} \subseteq G_n$} which involves all the variables \mbox{$x_1,\ldots,x_n$},
the following new system
\[
\left(\bigcup_{x_i \cdot x_j=x_k \in S} \{x_i \cdot x_j=x_k\}\right) \cup
\left\{2^{\textstyle 2^{\textstyle x_k}}=y_k:~k \in \{1,\ldots,n\}\right\} \cup \bigcup_{x_i+1=x_k \in S} \{y_i \cdot y_i=y_k\}
\]
is equivalent to ${\cal S}$. If the system ${\cal S}$ has only finitely many solutions in \mbox{non-negative} integers \mbox{$x_1,\ldots,x_n$},
then the new system has only finitely many solutions in \mbox{non-negative} integers \mbox{$x_1,\ldots,x_n,y_1,\ldots,y_n$}.
\end{observation}
\begin{proof}
It follows from Lemma~\ref{lem1}.
\end{proof}
\begin{observation}\label{obs2}
For every positive integer $n$, the following system
\[
\left\{\begin{array}{rcl}
x_1 \cdot x_1 &=& x_1 \\
\forall i \in \{1,\ldots,n-1\} ~2^{\textstyle 2^{\textstyle x_i}} &=& x_{i+1} ~({\rm if~} n>1)
\end{array}\right.
\]
has exactly two solutions in \mbox{non-negative} integers,
namely \mbox{$(\theta(1),\ldots,\theta(n))$}
and \mbox{$(f(1),\ldots,f(n))$}. The second solution has greater height.
\end{observation}
\begin{conjecture}\label{con2}
If a system \mbox{${\cal S} \subseteq G_n$}
has only finitely many solutions in \mbox{non-negative} rationals \mbox{$x_1,\ldots,x_n$},
then each such solution \mbox{$(x_1,\ldots,x_n)$} satisfies \mbox{$h(x_1,\ldots,x_n) \leqslant f(2n)$}.
\end{conjecture}
\par
Observations~\ref{obs3} and \ref{obs4} justify Conjecture~\ref{con2}.
\begin{observation}\label{obs3}
For every system \mbox{${\cal S} \subseteq G_n$} which involves all the variables \mbox{$x_1,\ldots,x_n$},
the following new system
\[
{\cal S} \cup \left\{2^{\textstyle 2^{\textstyle [x_k]}}=y_k:~k \in \{1,\ldots,n\}\right\} \cup
\bigcup_{x_i+1=x_k \in {\cal S}} \left\{y_i \cdot y_i=y_k\right\}
\]
is equivalent to ${\cal S}$. If the system ${\cal S}$ has only finitely many solutions in \mbox{non-negative}
rationals \mbox{$x_1,\ldots,x_n$}, then the new system has only finitely many solutions in \mbox{non-negative}
rationals \mbox{$x_1,\ldots,x_n,y_1,\ldots,y_n$}.
\end{observation}
\begin{proof}
It follows from Lemma~\ref{lem2}.
\end{proof}
\begin{observation}\label{obs4}
For every positive integer $n$, the following system
\[
\left\{\begin{array}{rcl}
x_1 \cdot x_1 &=& x_1 \\
\forall i \in \{1,\ldots,n-1\} ~2^{\textstyle 2^{\textstyle [x_i]}} &=& x_{i+1} ~({\rm if~} n>1)
\end{array}\right.
\]
has exactly two solutions in \mbox{non-negative} rationals,
namely \mbox{$(\theta(1),\ldots,\theta(n))$}
and \mbox{$(f(1),\ldots,f(n))$}. The second solution has greater height.
\end{observation}
\section{Algebraic lemmas}
\begin{lemma}\label{lem6} (cf. \cite[p.~100]{Robinson})
For every \mbox{non-negative} real numbers \mbox{$x,y,z$}, \mbox{$x+y=z$} if and only if
\begin{equation}\label{equ1}
((z+1)x+1)((z+1)(y+1)+1)=(z+1)^2(x(y+1)+1)+1
\end{equation}
\end{lemma}
\begin{proof}
The left side of equation~(\ref{equ1}) minus the right side of equation~(\ref{equ1}) equals \mbox{$(z+1)(x+y-z)$}.
\end{proof}
\par
Let $\alpha$, $\beta$, and $\gamma$ denote variables.
\begin{lemma}\label{lem7}
In \mbox{non-negative} integers, the equation \mbox{$x+y=z$} is equivalent to a system
which consists of equations of the forms \mbox{$\alpha+1=\gamma$} and \mbox{$\alpha \cdot \beta=\gamma$}.
The same is true for \mbox{non-negative} rationals.
\end{lemma}
\begin{proof}
It follows from Lemma~\ref{lem6}.
\end{proof}
\begin{lemma}\label{lem8}
Let \mbox{$D(x_1,\ldots,x_p) \in {\Z}[x_1,\ldots,x_p]$}.
Assume that \mbox{${\rm deg}(D,x_i) \geqslant 1$} for each \mbox{$i \in \{1,\ldots,p\}$}. We can compute a positive
integer \mbox{$n>p$} and a system \mbox{${\cal T} \subseteq G_n$} which satisfies the following two conditions:
\vskip 0.2truecm
\noindent
{\tt Condition 1.} For every \mbox{non-negative} integers \mbox{$\tilde{x}_1,\ldots,\tilde{x}_p$},
\[
D(\tilde{x}_1,\ldots,\tilde{x}_p)=0 \Longleftrightarrow
\exists \tilde{x}_{p+1},\ldots,\tilde{x}_n \in \N ~(\tilde{x}_1,\ldots,\tilde{x}_p,\tilde{x}_{p+1},\ldots,\tilde{x}_n) ~solves~ {\cal T}
\]
{\tt Condition 2.} If \mbox{non-negative} integers \mbox{$\tilde{x}_1,\ldots,\tilde{x}_p$} satisfy
\mbox{$D(\tilde{x}_1,\ldots,\tilde{x}_p)=0$}, then there exists a unique tuple
\mbox{$(\tilde{x}_{p+1},\ldots,\tilde{x}_n) \in {\N}^{n-p}$} such that the tuple
\mbox{$(\tilde{x}_1,\ldots,\tilde{x}_p,\tilde{x}_{p+1},\ldots,\tilde{x}_n)$} solves ${\cal T}$.
\vskip 0.2truecm
\noindent
Conditions 1 and 2 imply that the equation \mbox{$D(x_1,\ldots,x_p)=0$} and the system ${\cal T}$ have
the same number of solutions in \mbox{non-negative} integers.
\end{lemma}
\begin{proof}
We write down the polynomial \mbox{$D(x_1,\ldots,x_p)$} and replace each coefficient by the successor
of its absolute value. Let \mbox{$\widetilde{D}(x_1,\ldots,x_p)$} denote the obtained polynomial.
The polynomials \mbox{$D(x_1,\ldots,x_p)+\widetilde{D}(x_1,\ldots,x_p)$} and \mbox{$\widetilde{D}(x_1,\ldots,x_p)$}
have positive integer coefficients. The equation \mbox{$D(x_1,\ldots,x_p)=0$} is equivalent to
\[
D(x_1,\ldots,x_p)+\widetilde{D}(x_1,\ldots,x_p)+1=\widetilde{D}(x_1,\ldots,x_p)+1
\]
There exist a positive integer $a$ and a finite \mbox{non-empty} list $A$ such that
\begin{equation}\label{equ2}
D(x_1,\ldots,x_p)+\widetilde{D}(x_1,\ldots,x_p)+1=
\Bigl(\Bigl(\Bigl(\sum_{\textstyle (i_1,j_1,\ldots,i_k,j_k) \in A}
x_{\textstyle i_1}^{\textstyle j_1} \cdot \ldots \cdot x_{\textstyle i_k}^{\textstyle j_k}\Bigr)+
\underbrace{1\Bigr)+\ldots\Bigr)+1}_{\textstyle a~{\rm units}}
\end{equation}
and all the numbers \mbox{$k,i_1,j_1,\ldots,i_k,j_k$} belong to \mbox{$\N \setminus \{0\}$}.
There exist a positive integer $b$ and a finite \mbox{non-empty} list $B$ such that
\begin{equation}\label{equ3}
\widetilde{D}(x_1,\ldots,x_p)+1=
\Bigl(\Bigl(\Bigl(\sum_{\textstyle (i_1,j_1,\ldots,i_k,j_k) \in B}
x_{\textstyle i_1}^{\textstyle j_1} \cdot \ldots \cdot x_{\textstyle i_k}^{\textstyle j_k}\Bigr)+
\underbrace{1\Bigr)+\ldots\Bigr)+1}_{\textstyle b~{\rm units}}
\end{equation}
and all the numbers \mbox{$k,i_1,j_1,\ldots,i_k,j_k$} belong to \mbox{$\N \setminus \{0\}$}.
By Lemma~\ref{lem7}, we can equivalently express the equality of the right
sides of equations (\ref{equ2}) and (\ref{equ3}) using only equations of the forms
\mbox{$\alpha+1=\gamma$} and \mbox{$\alpha \cdot \beta=\gamma$}.
Consequently, we can effectively find the system ${\cal T}$.
\end{proof}
\par
The next lemma is a rational analogue of Lemma~\ref{lem8}.
\begin{lemma}\label{lem9}
Let \mbox{$D(x_1,\ldots,x_p) \in {\Z}[x_1,\ldots,x_p]$}.
Assume that \mbox{${\rm deg}(D,x_i) \geqslant 1$} for each \mbox{$i \in \{1,\ldots,p\}$}. We can compute a positive
integer \mbox{$n>p$} and a system \mbox{${\cal T} \subseteq G_n$} which satisfies the following two conditions:
\vskip 0.2truecm
\noindent
{\tt Condition 3.} For every \mbox{non-negative} rationals \mbox{$\tilde{x}_1,\ldots,\tilde{x}_p$},
\[
D(\tilde{x}_1,\ldots,\tilde{x}_p)=0 \Longleftrightarrow
\exists \tilde{x}_{p+1},\ldots,\tilde{x}_n \in {\Q} \cap [0,\infty) ~(\tilde{x}_1,\ldots,\tilde{x}_p,\tilde{x}_{p+1},\ldots,\tilde{x}_n) ~solves~ {\cal T}
\]
{\tt Condition 4.} If \mbox{non-negative} rationals \mbox{$\tilde{x}_1,\ldots,\tilde{x}_p$} satisfy
\mbox{$D(\tilde{x}_1,\ldots,\tilde{x}_p)=0$}, then there exists a unique tuple
\mbox{$(\tilde{x}_{p+1},\ldots,\tilde{x}_n) \in {\left(\Q \cap [0,\infty)\right)}^{n-p}$} such that the tuple
\mbox{$(\tilde{x}_1,\ldots,\tilde{x}_p,\tilde{x}_{p+1},\ldots,\tilde{x}_n)$} solves~${\cal T}$.
\vskip 0.2truecm
\noindent
Conditions 3 and 4 imply that the equation \mbox{$D(x_1,\ldots,x_p)=0$} and the system ${\cal T}$ have
the same number of solutions in \mbox{non-negative} rationals.
\end{lemma}
\section{Conjectural upper bounds on the heights of the solutions}
\begin{theorem}\label{the1}
If we assume Conjecture~\ref{con1} and a Diophantine equation \mbox{$D(x_1,\ldots,x_p)=0$} has only
finitely many solutions in \mbox{non-negative} integers, then an upper bound for these solutions
can be computed.
\end{theorem}
\begin{proof}
It follows from Lemma~\ref{lem8}.
\end{proof}
\par
The authors of \cite{DMR} proposed the following conjecture~\mbox{(\cite[p.~372]{DMR})}:
{\em there is no algorithm for listing the Diophantine equations
with infinitely many solutions in $\N$}. By Theorem~\ref{the1}, Conjecture~\ref{con1} guarantees
that such an algorithm exists.
\begin{theorem}\label{the3}
If we assume Conjecture~\ref{con1} and a Diophantine equation \mbox{$D(x_1,\ldots,x_p)=0$}
has only finitely many solutions in positive integers, then an upper bound for
these solutions can be computed.
\end{theorem}
\begin{proof}
We apply Theorem~\ref{the1} to the equation \mbox{$D(x_{1}+1,\ldots,x_{p}+1)=0$}.
Next, we increase the computed bound by $1$.
\end{proof}
\begin{theorem}\label{the4}
If we assume Conjecture~\ref{con1} and a Diophantine equation \mbox{$D(x_1,\ldots,x_p)=0$} has only finitely many
integer solutions, then an upper bound for their moduli can be computed by applying Theorem~\ref{the1}
to the equation
\[
\prod\limits_{\textstyle (i_1,\ldots,i_p) \in \{1,2\}^p}
D((-1)^{\textstyle i_1} \cdot x_1,\ldots,(-1)^{\textstyle i_p} \cdot x_p)=0
\]
\end{theorem}
\begin{lemma}\label{lem10} (\cite[Corollary~2,~p.~25]{Sierpinski})
If $a$ and $b$ are two relatively prime positive integers, then every integer \mbox{$n>ab$} can
be written in the form \mbox{$n=ax+by$}, where $x$, $y$ are positive integers.
\end{lemma}
\begin{lemma}\label{lem11}
For every \mbox{non-negative} integers $c$ and $d$, the following system
\begin{equation}\label{equ4}
\left\{\begin{array}{rcl}
cx+(d+1)y &=& (d+1)c+1 \\
x+y+u &=& (d+1)c+1
\end{array}\right.
\end{equation}
has at most finitely many solutions in \mbox{non-negative} integers $x$, $y$, $u$.
For every \mbox{non-negative} integers $c$ and $d$, system~(\ref{equ4}) is solvable in
\mbox{non-negative} integers $x$, $y$, $u$ if and only if $c$ and \mbox{$d+1$} are relatively prime.
\end{lemma}
\begin{proof}
The equality \mbox{$x+y+u=(d+1)c+1$} implies that \mbox{$x,y,u \leqslant (d+1)c+1$}.
Hence, at most finitely many \mbox{non-negative} integers \mbox{$x,y,u$} satisfy system~(\ref{equ4}).
The equality
\begin{equation}\label{equ5}
cx+(d+1)y=(d+1)c+1 
\end{equation}
gives \mbox{$cx+(d+1)(y-c)=1$}. Hence, the integers $c$ and \mbox{$d+1$} are relatively prime.
Conversely, assume that $c$ and \mbox{$d+1$} are relatively prime. By this, if \mbox{$c=0$},
then \mbox{$d=0$}. In this case, system~(\ref{equ4}) has exactly one solution in \mbox{non-negative} integers, namely
\mbox{$\left\{\begin{array}{rcl}
x &=& 0 \\
y &=& 1 \\
u &=& 0
\end{array}\right.$}.
If \mbox{$c>0$}, then Lemma~\ref{lem10} implies that there exist positive integers $x$ and $y$ that
satisfy equation~(\ref{equ5}). We set \mbox{$u=(c-1)x+dy$}. Then,
\[
x+y+u=x+y+(c-1)x+dy=cx+(d+1)y=(d+1)c+1
\]
\end{proof}
\begin{theorem}\label{the6}
Conjecture~\ref{con1} implies that there exists a computable upper bound
on the heights of the rationals that solve a Diophantine equation
with a finite number of solutions.
\end{theorem}
\begin{proof}
Let \mbox{$W(x_1,\ldots,x_n) \in {\Z}[x_1,\ldots,x_n]$}, and let
\[
\widehat{W}(x_1,\ldots,x_n)=\prod\limits_{\textstyle (i_1,\ldots,i_n) \in \{1,2\}^n}
W((-1)^{\textstyle i_1} \cdot x_1,\ldots,(-1)^{\textstyle i_n} \cdot x_n)
\]
If the equation \mbox{$W(x_1,\ldots,x_n)=0$} has only finitely many solutions in rationals \mbox{$x_1,\ldots,x_n$},
then the equation \mbox{$\widehat{W}(x_1,\ldots,x_n)=0$} has only finitely many solutions in \mbox{non-negative}
rationals \mbox{$x_1,\ldots,x_n$}. By Lemma~\ref{lem11}, it means that the system
\begin{equation}\label{equ6}
\left\{\begin{array}{rcl}
\widehat{W}\left(\frac{\textstyle y_1}{\textstyle z_1+1},\ldots,\frac{\textstyle y_n}{\textstyle z_n+1}\right) &=& 0 \\
\forall i \in \{1,\ldots,n\} ~~~~y_{i}s_{i}+(z_{i}+1)t_{i} &=& (z_{i}+1)y_{i}+1 \\
\forall i \in \{1,\ldots,n\} ~~~~s_{i}+t_{i}+u_{i} &=& (z_{i}+1)y_{i}+1
\end{array}\right.
\end{equation}
has only finitely many solutions in \mbox{non-negative} integers \mbox{$y_1,z_1,s_1,t_1,u_1,\ldots,y_n,z_n,s_n,t_n,u_n$}.
System~(\ref{equ6}) is equivalent to a single Diophantine equation.
By Lemma~\ref{lem8}, this equation is equivalent to a system of equations of the forms
\mbox{$\alpha \cdot \beta=\gamma$} and \mbox{$\alpha+1=\gamma$}. Next, we apply Theorem~\ref{the1}.
\end{proof}
\begin{theorem}\label{the2}
If we assume Conjecture~\ref{con2} and a Diophantine equation \mbox{$D(x_1,\ldots,x_p)=0$} has only
finitely many solutions in \mbox{non-negative} rationals, then an upper bound for their heights can be computed.
\end{theorem}
\begin{proof}
It follows from Lemma~\ref{lem9}.
\end{proof}
\begin{theorem}\label{the5}
If we assume Conjecture~\ref{con2} and a Diophantine equation \mbox{$D(x_1,\ldots,x_p)=0$} has only finitely many
rational solutions, then an upper bound for their heights can be computed by applying Theorem~\ref{the2}
to the equation
\[
\prod\limits_{\textstyle (i_1,\ldots,i_p) \in \{1,2\}^p}
D((-1)^{\textstyle i_1} \cdot x_1,\ldots,(-1)^{\textstyle i_p} \cdot x_p)=0
\]
\end{theorem}
\section{Finite-fold Diophantine representations}
The Davis-Putnam-Robinson-Matiyasevich theorem states that every recursively
enumerable set \mbox{${\cal M} \subseteq {\N}^n$} has a Diophantine
representation, that is
\[
(a_1,\ldots,a_n) \in {\cal M} \Longleftrightarrow
\exists x_1, \ldots, x_m \in \N ~~W(a_1,\ldots,a_n,x_1,\ldots,x_m)=0 \tag*{\texttt{(R)}}
\]
for some polynomial $W$ with integer coefficients, see \cite{Matiyasevich1}.
The polynomial $W$ can be computed, if we know the Turing \mbox{machine $M$} such
that, for all \mbox{$(a_1,\ldots,a_n) \in {\N}^n$}, $M$ halts on \mbox{$(a_1,\ldots,a_n)$} if
and only if \mbox{$(a_1,\ldots,a_n) \in {\cal M}$}, \mbox{see \cite{Matiyasevich1}}.
The representation~\texttt{(R)} is said to be \mbox{finite-fold}, if for every
\mbox{$a_1,\ldots,a_n \in \N$} the equation
\mbox{$W(a_1,\ldots,a_n,x_1,\ldots,x_m)=0$} has only finitely many solutions
\mbox{$(x_1,\ldots,x_m) \in {\N}^m$}. \mbox{Yu. Matiyasevich} conjectured that
each recursively enumerable set \mbox{${\cal M} \subseteq {\N}^n$} has a
\mbox{finite-fold} Diophantine representation, see \mbox{\cite[pp.~341--342]{DMR}},
\mbox{\cite[p.~42]{Matiyasevich2}}, and \mbox{\cite[p.~745]{Matiyasevich3}}.
Currently, he seems agnostic on his
 conjecture, see \mbox{\cite[p.~749]{Matiyasevich3}}.
In \mbox{\cite[p.~581]{Tyszka2}}, the author explains why Matiyasevich's conjecture although widely known is less widely accepted.
Matiyasevich's conjecture implies a negative answer to Open Problem~\ref{pro1}, see \mbox{\cite[p.~42]{Matiyasevich2}}.
\begin{lemma}\label{lem12}
Let \mbox{$W(x,x_1,\ldots,x_m) \in {\Z}[x,x_1,\ldots,x_m]$}. We claim that the function
\[
\N \ni b \mapsto W(b,x_1,\ldots,x_m) \in {\Z}[x_1,\ldots,x_m]
\]
is computable.
\end{lemma}
\begin{theorem}\label{the7}
Conjecture~\ref{con1} implies that if a set \mbox{${\cal M} \subseteq \N$} has a \mbox{finite-fold}
Diophantine representation, then ${\cal M}$ is computable.
\end{theorem}
\begin{proof}
Let a set \mbox{${\cal M} \subseteq \N$} have a \mbox{finite-fold} Diophantine representation.
It means that there exists a polynomial \mbox{$W(x,x_1,\ldots,x_m)$} with integer coefficients such that
\[
\forall b \in \N~ \Bigl(b \in {\cal M} \Longleftrightarrow \exists x_1, \ldots, x_m \in \N ~~W(b,x_1,\ldots,x_m)=0\Bigr)
\]
and for every \mbox{$b \in \N$} the equation \mbox{$W(b,x_1,\ldots,x_m)=0$} has only finitely many solutions
\mbox{$(x_1,\ldots,x_m) \in {\N}^m$}. By Lemma~\ref{lem12} and Theorem~\ref{the1}, there is a computable function \mbox{$\xi \colon \N \to \N$}
such that for each \mbox{$b,x_1,\ldots,x_m \in \N$} the equality \mbox{$W(b,x_1,\ldots,x_m)=0$} implies
\mbox{$\max(x_1,\ldots,x_m) \leqslant \xi(b)$}.
Hence, we can decide whether or not a \mbox{non-negative} integer $b$ belongs to ${\cal M}$ by checking
whether or not the equation \mbox{$W(b,x_1,\ldots,x_m)=0$} has an integer solution in the box \mbox{$[0,\xi(b)]^m$}.
\end{proof}
\par
Theorem~\ref{the7} remains true if we change the bound \mbox{$f(2n)$} in Conjecture~\ref{con1}
to any other computable bound \mbox{$\delta(n)$}.
\section{Theorems on relative decidability}
\begin{conjecture}\label{friedman1} (\cite{Friedman})
The set of all Diophantine equations which have only finitely many rational solutions is not recursively enumerable.
\end{conjecture}
\par
Conjecture~\ref{friedman1} implies Conjecture~\ref{friedman2}.
\begin{conjecture}\label{friedman2}
The set of all Diophantine equations which have only finitely many rational solutions is not computable.
\end{conjecture}
\par
By Theorem~\ref{the6}, Conjecture~\ref{con1} implies Conjecture~\ref{maj2}.
By Theorem~\ref{the5}, Conjecture~\ref{con2} implies Conjecture~\ref{maj2}.
\begin{conjecture}\label{maj2}
There exists an algorithm which takes as input a Diophantine equation \mbox{$D(x_1,\ldots,x_p)=0$} and
returns an integer \mbox{$b \geqslant 2$}, where $b$ is greater than the number of rational solutions,
if the solution set is finite.
\end{conjecture}
\vskip 0.01truecm
\noindent
{\bf Guess}~\mbox{(\cite[p.~16]{Kim})}. {\em The question whether or not a Diophantine
equation has only finitely many rational solutions is decidable with an oracle
for deciding whether or not a Diophantine equation has a rational solution.}
\vskip 0.2truecm
\par
Originally, Minhyong Kim formulated the Guess as follows: for rational solutions, the finiteness
problem is decidable relative to the existence problem.
Conjecture~\ref{friedman2} and the Guess imply that there is no algorithm
which decides whether or not a Diophantine equation has a rational solution.
\begin{theorem}\label{the8}
Conjecture~\ref{maj2} implies the Guess.
\end{theorem}
\begin{proof}
Assuming that Conjecture~\ref{maj2} holds,
the execution of Flowchart 1 decides whether or not a Diophantine equation
\mbox{$D(x_1,\ldots,x_p)=0$} has only finitely many rational solutions.
\vskip 0.2truecm
\begin{center}
\begin{tikzpicture}[very thick]
\ttfamily
\node at (6.5,10.3) {Start};
\node at (6.5,9.15) {Input a Diophantine equation $D \left( x_1, \ldots, x_p
\right) = 0$};
\node at (6.5,8) {Compute the bound $b$};
\node[text width=11cm, text centered] at (6.5,6.05) {Does the equation\break
$\displaystyle \left( \sum_{k=1}^{b} D^{2}\left( x_{1, k}, \ldots, x_{p, k}
\right) \right) + \left( \left( \prod_{1 \leqslant u < v \leqslant b}
\sum_{i=1}^{p} \left( x_{i, u} - x_{i, v} \right)^{2} \right) \cdot y - 1
\right)^{2}\!\!=\!\!0$\break have a rational solution?};
\node[text width=11cm, text centered] at (6.5,3.6) {Print "The equation $D\left(
x_1, \ldots, x_p \right) = 0$ has infinitely many rational solutions"};
\node[text width=11cm, text centered] at (6.5,1.9) {Print "The equation $D\left(
x_1, \ldots, x_p \right) = 0$ has only finitely many rational solutions"};
\node at (6.5,.3) {Stop};
\draw (6.5,10.3) ellipse(.7cm and 0.3cm);
\draw (1.3,8.8) -- (11.3,8.8) -- (11.7,9.5) -- (1.7,9.5) -- cycle;
\draw (4.3,7.7) rectangle (8.7,8.3);
\draw (1.1,4.9) rectangle (11.9,7.2);
\draw (1.9,3) -- (10.7,3) -- (11.1,4.2) -- (2.3,4.2) -- cycle;
\draw (1.9,1.3) -- (10.7,1.3) -- (11.1,2.5) -- (2.3,2.5) -- cycle;
\draw (6.5,.3) ellipse(.76cm and 0.3cm);
\draw[->] (6.5,10) -- (6.5,9.5);
\draw[->] (6.5,8.8) -- (6.5,8.3);
\draw[->] (6.5,7.7) -- (6.5,7.2);
\draw[->] (6.5,4.9) -- (6.5,4.2) node[midway,right] {Yes};
\draw[->] (6.5,1.3) -- (6.5,.6);
\draw[->] (11.9,6.05) -- (12.6,6.05) -- (12.6,1.9) -- (10.9,1.9);
\node[above right] at (11.9,6.05) {No};
\draw[->] (2.1,3.6) -- (1.1,3.6) -- (1.1,1) -- (6.5,1);
\end{tikzpicture}
\end{center}
\vskip 0.2truecm
\centerline{{\it Flowchart 1: Conjecture~\ref{maj2} implies the Guess}}
\end{proof}
\begin{lemma}\label{xyz}
A Diophantine equation \mbox{$D(x_1,\ldots,x_p)=0$} has no solutions in rationals \mbox{$x_1,\ldots,x_p$}
if and only if the equation \mbox{$D(x_1,\ldots,x_p)+0 \cdot x_{p+1}=0$} has only finitely many
solutions in rationals \mbox{$x_1,\ldots,x_p,x_{p+1}$}.
\end{lemma}
\begin{theorem}\label{the9}
If the set of all Diophantine equations which
have only finitely many rational solutions is recursively enumerable, then there exists an algorithm
which decides whether or not a Diophantine equation has a rational solution.
\end{theorem}
\begin{proof}
We find a computable surjective function \mbox{$\eta \colon \N \to \bigcup_{n=1}^\infty\limits {\Q}^n$}.
Suppose that \mbox{$\{{\cal A}_n=0\}_{n=0}^\infty$} is a computable sequence of all Diophantine
equations which have only finitely many rational solutions.
By Lemma~\ref{xyz}, the execution of Flowchart 2 decides whether or not a Diophantine
equation \mbox{$D(x_1,\ldots,x_p)=0$} has a rational solution.
\begin{center}
\includegraphics[scale=1.5]{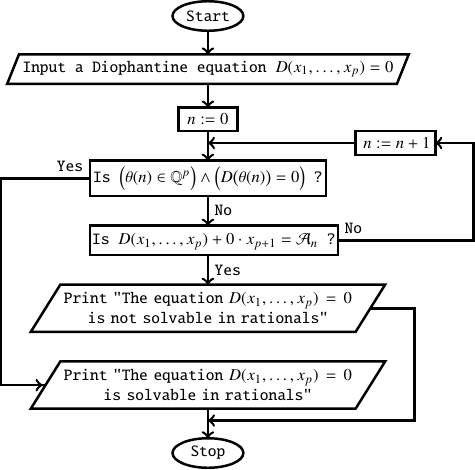}
\end{center}
\vskip 0.01truecm
\centerline{{\it Flowchart 2: An algorithm which conditionally decides whether or not}}
\vskip 0.01truecm
\centerline{{\it a Diophantine equation $D(x_1,\ldots,x_p)=0$ has a rational solution}}
\end{proof}
\par
By Theorem~\ref{the4}, Conjecture~\ref{con1} implies Conjecture~\ref{zwyc}.
\begin{conjecture}\label{zwyc}
There exists an algorithm which takes as input a Diophantine equation \mbox{$D(x_1,\ldots,x_p)=0$} and
returns an integer \mbox{$r \geqslant 2$},
where $r$ is greater than the moduli of integer solutions, if the solution set is finite.
\end{conjecture}
\par
Conjecture~\ref{zwyc} implies Conjecture~\ref{maj22}.
\begin{conjecture}\label{maj22}
There exists an algorithm which takes as input a Diophantine equation \mbox{$D(x_1,\ldots,x_p)=0$} and
returns an integer \mbox{$c \geqslant 2$}, where $c$ is greater than the number of integer solutions,
if the solution set is finite.
\end{conjecture}
\begin{theorem}\label{t9}
Conjecture~\ref{maj22} implies Conjecture~\ref{zwyc}.
\end{theorem}
\begin{proof}
Assume that a Diophantine equation \mbox{$D(x_1,\ldots,x_p)=0$} has only finitely many integer solutions.
Then, the equation
\begin{equation}\label{seven}
D^2\left(x_1,\ldots,x_p\right)+\left(\left(x_1^2+\ldots+x_p^2\right)-\left(y_1^2+y_2^2+y_3^2+y_4^2\right)-\left(z_1^2+z_2^2+z_3^2+z_4^2\right)\right)^2=0
\end{equation}
has only finitely many integer solutions. By Lagrange's \mbox{four-square} theorem,
every integer \mbox{$p$-tuple} \mbox{$(a_1,\ldots,a_p)$} with \mbox{$D(a_1,\ldots,a_p)=0$}
implies the existence of \mbox{$a_1^2+\ldots+a_p^2+1$} integer solutions of equation~(\ref{seven}).
The inequality \mbox{$|a_1|,\ldots,|a_p|<a_1^2+\ldots+a_p^2+1$} completes the proof.
\end{proof}
\newpage
\begin{theorem}\label{the11}
Conjecture~\ref{maj22} implies that the question whether or not a Diophantine equation has only
finitely many integer solutions is decidable with an oracle for deciding whether or not a Diophantine equation has
an integer solution.
\end{theorem}
\begin{proof}
Assuming that Conjecture~\ref{maj22} holds, Lagrange's \mbox{four-square} theorem guarantees that the
execution of Flowchart 3 decides whether or not a Diophantine equation
\mbox{$D(x_1,\ldots,x_p)=0$} has only finitely many integer solutions.
\vskip 0.2truecm
\begin{center}
\begin{tikzpicture}[very thick]
\ttfamily
\node at (7.5,10.3) {Start};
\node at (7.5,9.15) {Input a Diophantine equation $D \left( x_1, \ldots, x_p
\right) = 0$};
\node at (7.5,8) {Compute the bound $c$};
\node[text width=14cm, text centered] at (7.45,6.05) {Does the equation\break
$\displaystyle \left( \sum_{k=1}^c D^{2}\left( x_{1, k}, \ldots,
x_{p, k} \right) \right) + \left( \left( \prod_{1 \leqslant u < v \leqslant c}
\sum_{i=1}^{p} \left( x_{i, u} - x_{i, v} \right)^{2} \right)
- s^{2} - t^{2} - u^{2} - v^{2} - 1 \right)^{2}\!\!=\!\!0$\break have an integer
solution?};
\node[text width=11cm, text centered] at (7.5,3.6) {Print "The equation $D\left(
x_1, \ldots, x_p \right) = 0$ has infinitely many integer solutions"};
\node[text width=11cm, text centered] at (7.5,1.9) {Print "The equation $D\left(
x_1, \ldots, x_p \right) = 0$ has only finitely many integer solutions"};
\node at (7.5,.3) {Stop};
\draw (7.5,10.3) ellipse(.7cm and 0.3cm);
\draw (2.3,8.8) -- (12.3,8.8) -- (12.7,9.5) -- (2.7,9.5) -- cycle;
\draw (5.3,7.7) rectangle (9.7,8.3);
\draw (.7,4.9) rectangle (14.2,7.2);
\draw (2.9,3) -- (11.7,3) -- (12.1,4.2) -- (3.3,4.2) -- cycle;
\draw (2.9,1.3) -- (11.7,1.3) -- (12.1,2.5) -- (3.3,2.5) -- cycle;
\draw (7.5,.3) ellipse(.76cm and 0.3cm);
\draw[->] (7.5,10) -- (7.5,9.5);
\draw[->] (7.5,8.8) -- (7.5,8.3);
\draw[->] (7.5,7.7) -- (7.5,7.2);
\draw[->] (7.5,4.9) -- (7.5,4.2) node[midway,right] {Yes};
\draw[->] (7.5,1.3) -- (7.5,.6);
\draw[->] (14.2,6.05) -- (14.8,6.05) -- (14.8,1.9) -- (11.9,1.9);
\node[above right] at (14.2,6.05) {No};
\draw[->] (3.1,3.6) -- (2.1,3.6) -- (2.1,1) -- (7.5,1);
\end{tikzpicture}
\vskip 0.2truecm
\centerline{{\it Flowchart 3: Conjecture~\ref{maj22} implies the Guess reformulated for integer solutions}}
\end{center}
\end{proof}
\par
The conclusion of Theorem~\ref{the11} contradicts the following conjecture of Yuri Matiyasevich~\mbox{(\cite[p.~16]{Kim})}:
{\em the finiteness problem for integral points is undecidable relative to the existence problem}.
\section{A summary of the main theorems}
The main results of this article are summarized in Flowchart 4.
\begin{center}
\begin{tikzpicture}[ultra thick]
\tt
\node[text width=6.5cm, text centered] at (3.35,21.83) {Yuri Matiyasevich's
conjecture on finite-fold Diophantine representation is false};
\node[text width=8.65cm, text centered] at (11.35,20.6) {Only computable subsets
of $\mathbb{N}$ may have a~finite-fold Diophantine representation};
\node[text width=5.9cm, text centered] at (3.7,18.9) {There exists an algorithm for
listing the Diophantine equations with infinitely many solutions in
$\mathbb{N}$};
\node[text width=5.4cm, text centered] at (3.45,16.45) {There exists a
computable upper bound on the number of integer solutions};
\node[text width=6cm, text centered] at (11.5,16.45) {There exists a computable
upper bound on the moduli of integer solutions};
\node[text width=5.5cm, text centered] at (11.45,13.8) {There exists a computable
upper bound on the heights of rational solutions};
\node[text width=8.0cm, text centered] at (11.75,11.65) {There exists a
computable upper bound on the number of rational solutions};
\node[text width=8.0cm, text centered] at (4.04,8) {Minhyong Kim's Guess
reformulated for integer solutions:\\ \em{the question whether or not
a~Diophantine equation has only finitely many integer solutions is decidable
with an oracle for deciding whether or not a~Diophantine equation has an integer
solution}};
\node[text width=6.35cm, text centered] at (12.47,7.77) {\hspace{6mm}Minhyong
Kim's Guess\\ ([2, p. 16]):\\ \em{the question whether or not a Diophantine
equation has only finitely many rational solutions is decidable with an~oracle
for deciding whether or not a Diophantine equation has a rational solution}};
\node[text width=7.6cm, text centered] at (3.8,4.05) {There is no algorithm
which decides whether or not a Diophantine equation has a rational solution};
\node[text width=4.55cm, text centered] at (13.4,4.05) {The conjunction of the
statements A and B};
\node[text width=8.4cm, text centered] at (4.3,1.1) {Harvey Friedman's
conjecture:\\ \em{the set of all Diophantine equations which have only finitely
many rational solutions is not recursively enumerable}};
\node[text width=5.4cm, text centered] at (13.1,1.35) {\hspace{8mm}The set of
all\\\hspace*{5mm}Diophantine equations which have only finitely many rational
solutions
is not computable};
\node at (11.5,18.65) {Conjecture~\ref{con1}};
\node at (4.75,13.8) {Conjecture~\ref{con2}};
\draw (0,21) rectangle (6.7,22.7);
\draw (6.95,20) rectangle (15.8,21.2);
\draw (.7,17.8) rectangle (6.7,20);
\draw (10.05,18.3) rectangle (12.95,19);
\draw (.7,15.6) rectangle (6.2,17.3);
\draw (8.7,15.6) rectangle (14.3,17.3);
\draw (3.3,13.45) rectangle (6.2,14.15);
\draw (8.6,13) rectangle (14.3,14.6);
\draw (7.7,11.1) rectangle (15.8,12.2);
\draw (0,5.9) rectangle (8.1,10.1);
\draw (9.2,5.4) rectangle (15.8,10.1);
\draw (0,3.2) rectangle (7.6,4.9);
\draw (11,3.5) rectangle (15.8,4.6);
\draw (0,0) rectangle (8.6,2.2);
\draw (10.4,0) rectangle (15.8,2.7);
\draw[->] (2.3,15.6) -- (2.3,10.1);
\draw[->] (8.7,16.2) -- (6.2,16.2);
\draw[->] (11,4.05) -- (7.6,4.05);
\draw[->] (8.6,1.1) -- (10.4,1.1);
\draw[->] (12,13) -- (12,12.2);
\draw[->] (12,4.6) -- (12,5.4);
\draw[->] (12,3.5) -- (12,2.7);
\draw[->] (.5,10.1) -- (.5,21);
\draw[->] (2.3,3.2) -- (2.3,2.2) node[midway, right] {Theorem~\ref{the9}};
\draw[->] (7.6,3.2) -- (10.4,2.7) node[midway,rotate=-12,above] {by Lemma 10};
\draw[->] (2.3,20) -- (2.3,21) node[midway, right] {\mbox{\cite[p.~372]{DMR}}};
\draw[->] (12,15.6) -- (12,14.6) node[midway, right] {Lemma 3 in \cite{Tyszka0}};
\draw[->] (10.05,18.65) -- (6.7,18.65) node[midway, above] {by Theorem~\ref{the1}};
\draw[->] (12,19) -- (12,20) node[midway, right] {Theorem~\ref{the7}};
\draw[->] (12,18.3) -- (12,17.3) node[midway, right] {Theorem~\ref{the4}};
\draw[->] (12,11.1) -- (12,10.1) node[midway, right] {Theorem~\ref{the8}};
\draw[->] (6.2,16.7) -- (8.7,16.7) node[midway, above] {Theorem~\ref{t9}};
\draw[->] (6.2,13.8) -- (8.6,13.8) node[midway, above] {Theorem~\ref{the5}};
\node[right] at (2.3,11.65) {Theorem~\ref{the11}};
\node[right] at (.5,12.85) {\rotatebox{90}{\mbox{\cite[p.~372]{DMR}}}};
\draw (12.95,18.65) -- (15.3,18.65) -- (15.3,14.05);
\draw (15.3,18.45) -- (15.3,14.15) node[midway, right] {\rotatebox{270}{Theorem~\ref{the6}}};
\draw[->] (15.3,14.35) -- (15.3,13.95) -- (14.3,13.95);
\draw[->] (12,21.2) -- (12,22.2) -- (6.7,22.2);
\node[below right] at (9.2,10.1) {\Large{A}};
\draw (9.2,8.9) -- (10.4,10.1);
\node[below right] at (10.4,2.7) {\Large{B}};
\draw (10.4,1.5) -- (11.6,2.7);
\end{tikzpicture}
\end{center}
\vskip 0.1truecm
\centerline{{\it Flowchart 4: A flowchart of implications which illustrates the main theorems}}

\noindent
Apoloniusz Tyszka\\
University of Agriculture\\
Faculty of Production and Power Engineering\\
Balicka 116B, 30-149 Krak\'ow, Poland\\
E-mail: \url{rttyszka@cyf-kr.edu.pl}
\end{sloppypar}

\begin{thebibliography}{11}

\bibitem{DMR}
\mbox{M. Davis}, \mbox{Yu. Matiyasevich}, \mbox{J. Robinson},
\newblock
{\em Hilbert's tenth problem. Diophantine equations: positive aspects of a negative solution;}
in: Mathematical developments arising from Hilbert problems (ed.~F.~E.~Browder),
\newblock
Proc. Sympos. Pure Math., vol. 28, Part 2, Amer. Math. Soc., Providence, RI, 1976, \mbox{323--378},
\newblock
\url{http://dx.doi.org/10.1090/pspum/028.2};
\newblock
reprinted in: The collected works of Julia Robinson (ed.~S.~Feferman), Amer. Math. Soc.,
Providence, RI, 1996, \mbox{269--324}.

\bibitem{Friedman}
\mbox{H. Friedman},
\newblock
{\em Complexity of statements,}
\newblock
April 20, 1998,
\newblock
\url{http://www.cs.nyu.edu/pipermail/fom/1998-April/001843.html}.

\bibitem{Kim}
\mbox{M. Kim},
\newblock
{\em On relative computability for curves,}
\newblock
Asia Pac. Math. Newsl. 3 (2013), \mbox{no. 2}, \mbox{16--20},
\newblock
\url{http://www.asiapacific-mathnews.com/03/0302/0016_0020.pdf}.

\bibitem{Matiyasevich1}
\mbox{Yu. Matiyasevich},
\newblock
{\em Hilbert's tenth problem,}
\newblock
MIT Press, Cambridge, MA, 1993.

\bibitem{Matiyasevich2}
\mbox{Yu. Matiyasevich},
\newblock
{\em Hilbert's tenth problem: what was done and what is to be done;}
\newblock
in: Hilbert's tenth problem: relations with arithmetic and algebraic geometry (Ghent, 1999),
\newblock
Contemp. Math. 270, \mbox{1--47}, Amer. Math. Soc., Providence, RI, 2000,
\newblock
\url{http://dx.doi.org/10.1090/conm/270}.

\bibitem{Matiyasevich3}
\mbox{Yu. Matiyasevich},
\newblock
{\em Towards \mbox{finite-fold} Diophantine representations,}
\newblock
J. Math. Sci. (N. Y.) \mbox{vol. 171}, \mbox{no. 6}, 2010, \mbox{745--752},\\
\newblock
\url{http://dx.doi.org/10.1007%2Fs10958-010-0179-4}.

\bibitem{Robinson}
\mbox{J. Robinson},
\newblock
{\em Definability and decision problems in arithmetic,}
\newblock
J. Symbolic Logic 14 (1949), \mbox{no.~2}, \mbox{98--114},
\url{http://dx.doi.org/10.2307/2266510};
\newblock
reprinted in: The collected works of Julia Robinson (ed.~S.~Feferman),
Amer. Math. Soc., Providence, RI, 1996, \mbox{7--23}.

\bibitem{Sierpinski}
\mbox{W. Sierpi{\'n}ski},
\newblock
{\em Elementary theory of numbers,} 2nd~ed. (ed.~A.~Schinzel),
\newblock
PWN (Polish Scientific Publishers) and North-Holland, Warsaw-Amsterdam, 1987.

\bibitem{Tyszka0}
\mbox{A. Tyszka},
\newblock
{\em Conjecturally computable functions which unconditionally
do not have any \mbox{finite-fold} Diophantine representation,}
\newblock
Inform. Process. Lett. 113 (2013), \mbox{no. 19--21}, \mbox{719--722},
\newblock
\url{http://dx.doi.org/10.1016/j.ipl.2013.07.004}.

\bibitem{Tyszka1}
\mbox{A. Tyszka},
\newblock
{\em A hypothetical way to compute an upper bound for the heights
of solutions of a Diophantine equation with a finite number of solutions,}
Proceedings of the 2015 Federated Conference on Computer Science and Information Systems
(eds. M. Ganzha, L. Maciaszek, M. Paprzycki);
{\em Annals of Computer Science and Information Systems}, \mbox{vol. 5}, \mbox{709--716},
IEEE Computer Society Press, 2015, \url{http://dx.doi.org/10.15439/2015F41}.

\bibitem{Tyszka2}
\mbox{A. Tyszka},
\newblock {\em All functions \mbox{$g \colon \N \to \N$} which have a \mbox{single-fold} Diophantine
representation are dominated by a \mbox{limit-computable} function \mbox{$f \colon \N \setminus \{0\} \to \N$}
which is implemented in {\sl MuPAD} and whose computability is an open problem;}
in: Computation, cryptography, and network security (eds. N.~J.~Daras and M.~Th.~Rassias),
\newblock
Springer, Cham, Switzerland, 2015, \mbox{577--590},
\newblock
\url{http://dx.doi.org/10.1007/978-3-319-18275-9_24}.

\end{thebibliography}
\end{document}